\documentclass[11pt]{article}
\usepackage[margin=1.3cm]{geometry}
\usepackage[T1]{fontenc}
\usepackage{natbib}
\usepackage[figuresright]{rotating}
\usepackage{float}
\usepackage{graphicx}
\usepackage{epstopdf}
\usepackage{url}
\usepackage[colorlinks,citecolor=blue,linkcolor=blue,urlcolor=red,pagebackref]{hyperref}
\usepackage{amsmath,amsthm}
\usepackage{longtable}
\usepackage[ruled,vlined]{algorithm2e}
\usepackage[dvipsnames]{xcolor}
\usepackage{orcidlink}

\theoremstyle{definition}

\theoremstyle{remark}

\providecommand{\keywords}[1]
{
  \small	
  \textbf{Keywords:} #1
}

\begin{document}

\title {Exploring the Performance of Genetic Algorithm and Variable Neighborhood Search for Solving the Single Depot Multiple Set Orienteering Problem: A Comparative Study}                      

\author{Ravi Kant\orcidlink{0000-0001-5835-4731}\footnote{Department of Computer Science and Information Systems, Birla Institute of Technology and Science, Pilani,
Pilani-333031, Rajasthan, I\textsc{ndia}. Email: \texttt{p20190020@pilani.bits-pilani.ac.in}.} \and 
Sarthak Agarwal\footnote{Department of Computer Science and Information Systems,
Birla Institute of Technology and Science, Pilani, Pilani-333031, Rajasthan, I\textsc{ndia}. Email: \texttt{asarthak2002@gmail.com}.} \and 
Aakash Gupta \footnote{Department of Computer Science and Information Systems,
Birla Institute of Technology and Science, Pilani, Pilani-333031, Rajasthan, I\textsc{ndia}. Email: \texttt{uchanahome8@gmail.com}.} \and 
Abhishek Mishra\orcidlink{0000-0002-2205-0514}\footnote{Department of Computer Science and Information Systems, Birla Institute of Technology and Science, Pilani, Pilani-333031, Rajasthan, I\textsc{ndia}. Email: \texttt{abhishek.mishra@pilani.bits-pilani.ac.in}.}}

\maketitle

\begin{abstract}

This article discusses the single Depot multiple Set Orienteering Problem (sDmSOP), a recently suggested generalization of the Set Orienteering Problem (SOP). This problem aims to discover a path for each traveler over a subset of vertices, where each vertex is associated with only one cluster, and the total profit made from the clusters visited is maximized while still fitting within the available budget constraints. The profit can be collected only by visiting at least one cluster vertex. According to the SOP, each vertex cluster must have at least one of its visits counted towards the profit for that cluster. Like to the SOP, the sDmSOP restricts the number of clusters visited based on the budget for tour expenses. To address this problem, we employ the Genetic Algorithm (GA) and Variable Neighborhood Search (VNS) meta-heuristic. The optimal solution for small-sized problems is also suggested by solving the Integer Linear Programming (ILP) formulation using the General Algebraic Modeling System (GAMS) 37.1.0 with CPLEX for the sDmSOP. Promising computational results are presented that demonstrate the practicability of the proposed GA, VNS meta-heuristic, and ILP formulation by demonstrating substantial improvements to the solutions generated by VNS than GA while simultaneously needing much less time to compute than CPLEX.

\end{abstract}


\maketitle

\section{Introduction}
\label{section:1}

The single Depot multiple Set Orienteering Problem (sDmSOP) proposed by \cite{kant2023single} extends the Set Orienteering Problem (SOP). It aims to identify multiple tours, with precisely one tour designated for each traveler, traversing a subset of clusters. These tours must start and end at the same depot (i.e., depot $1$ in case of the sDmSOP) while maximizing the collected profit within a predefined budget $B$ (which may vary depending on the problem, i.e., time or distance budget). Unlike the traditional SOP, which focuses on finding a single path for one traveler, the sDmSOP necessitates discovering multiple closed paths for each traveler. The sDmSOP is a generalization of the SOP in the same way as the multiple Traveling Salesman Problem (mTSP) is a generalization of the Traveling Salesman Problem (TSP). 

Just as the SOP finds numerous real-life applications,  the sDmSOP also proves invaluable in contexts like mass distribution products. Here, multiple carriers (travelers) opt to serve individual retailers within the supply chain (customer) clusters. Importantly, each supply chain is served only once, by the same carrier or different carriers. 
 
The practical utility of the sDmSOP extends beyond mass distribution scenarios. For example, in the domain of sales and distribution, multiple carriers may serve individual retailers within supply chain clusters. Once again, the critical constraint remains that each supply chain is serviced only once, whether by the same carrier or different carriers. Subsequently, products are internally distributed among all the retailers within the chain.

Based on our understanding, the SOP has been exclusively investigated in the study conducted by \cite{archetti2018set}, who also provided a comprehensive set of benchmark test instances containing up to $1084$ vertices. In their work, they proposed a mathematical formulation and a MAtheuristic to solve the SOP, called the MA-SOP. This meta-heuristic integrates a Mixed Integer Linear Programming (MILP) approach into a tabu-search framework. The MILP-based move is employed when the tabu-search procedure fails to identify a new non-tabu feasible solution. This entails removing a set of visited clusters from the tour and subsequently inserting a selection of non-visited clusters by solving the MILP reallocation model. After that, \cite{pvenivcka2019variable} introduced an innovative Integer Linear Programming (ILP) formulation for the SOP, requiring fewer decision variables compared to the mathematical formulation proposed by \cite{archetti2018set}. This formulation enables optimal solutions to be achieved within shorter computational times, particularly for small benchmark instances with (\( n \leq 198 \)) using CPLEX 12.6.1., here $n$ is the number of nodes in the instance. However, despite its effectiveness, this model still struggles to provide optimal solutions for instances where (\( n \geq 200 \)) within a 9-hour time limit except in a few instances. Consequently, the authors introduced a Variable Neighbourhood Search (VNS) meta-heuristic for the SOP. The VNS uses two local search operators (cluster move and cluster exchange) and two shake operators (path move and path exchange) to find the near-optimal solution of the SOP.

In recent developments, \cite{carrabs2021biased} proposed a Biased Random-Key Genetic Algorithm (BRKGA) for the SOP. This algorithm incorporates three local search procedures for offspring improvement alongside the integration of a hash table and a three-dimensional matrix to circumvent redundant computations. Additionally, \cite{golmankhaneh2020acs} devised an ant colony system meta-heuristic for the SOP, called ACS-SOP. \cite{dontas2023adaptive} introduced a novel Adaptive Memory Meta-Heuristic named AMMH for the SOP. This local search algorithm integrates mathematical programming components to tackle various sub-problems alongside an adaptive memory structure to generate high-quality initial solutions for the SOP. Finally, \cite{lu2024effective} proposed a Hybrid Evolutionary Algorithm (HEA) to address the SOP, merging an effective local refinement process with a cluster-centric crossover operator and a randomized mutation procedure. The local refinement alternates between feasible and infeasible local searches to enhance search intensification. Computational findings demonstrate a substantial superiority of this algorithm over the top-performing heuristics for the SOP. \cite{archetti2024new} proposed a new mathematical formulation for the SOP and a branch-and-cut algorithm to provide an exact solution for the problem. The branch-and-cut approach effectively solves instances with up to 200 vertices. The algorithm cannot solve the SOP instances with more than 107 nodes when the value of $w$ is higher (i.e., larger budget).

In summary, further research on practical algorithms for the SOP remains imperative. While current mathematical models can only solve small instances, existing heuristics have struggled to perform well on larger instances consistently. For instance, while BRKGA has shown promise as the fastest heuristic for certain large SOP instances, VNS-SOP offers significantly faster computation at the cost of solution quality. MA-SOP and AMMH are time-consuming due to their combined use of mathematical programming and local search, while ACS-SOP presents even greater computational complexity.

This study presents the Genetic Algorithm (GA) and VNS meta-heuristic for addressing the sDmSOP. The effectiveness of these algorithms is evaluated using instances of the Generalized Traveling Salesman Problem (GTSP) featuring up to $1084$ vertices sourced from benchmark examples. The study compares the algorithmic outcomes against optimal solutions, showcasing their performance in more minor instances. The validation of optimal solutions for these smaller instances is carried out through the General Algebraic Modeling System (GAMS) and CPLEX.

The rest of the paper is structured as follows: Section \ref{section:2} presents a comprehensive problem description and mathematical formulation of the sDmSOP. Following this, Sections \ref{section:3} and \ref{section:4} investigate an in-depth exploration of the GA and VNS. Moving forward to Section \ref{section:5}, computational tests are conducted, and comparisons with the results obtained from GAMS, GA, and VNS are carefully examined. Section \ref{section:6} includes the conclusions drawn from the study. 

\section{Problem Description and Formulation}
\label{section:2}

The sDmSOP is a generalization of the SOP. Hence, the SOP must be formally defined first. The SOP can be formalized on a directed complete graph $G (V, E)$, where $V=\{\, 1, \ldots, n \,\}$ is the set of vertices and $E= \left \{\, (i, j) \mid (i, j) \in V^2 \, \right \}$ is the set of edges. The edge $(i, j)$ is defined as an edge from the vertex $i$ to the vertex $j$; moreover, a cost $c_{ij}\geq 0$, is associated with the edge $(i, j)$. The vertices are partitioned into disjoint sets $S=\{\, s_{1}, \ldots, s_{p} \,\}$ such that their union contains all the vertices of the graph. The objective of the problem is to gain maximum profit by visiting the possible number of sets within a distance constraint $B$ with the predefined starting and ending depot. The profit from a set can be collected if only one vertex of a set is visited by a traveler.

The sDmSOP can be formalized on a weighted directed complete graph $G (V, E)$ with weights $\{\, c_{ij} \mid (i, j) \in E \,\}$ as defined above. The vertices are partitioned into disjoint sets $S=\{\, s_{1},\cdots, s_{p} \,\}$ such that their union contains all vertices of the graph. The objective of the problem is to gain maximum profit by visiting the possible number of sets within a distance constraint ($B$) using $t \in \{\, 1,\ldots, m \,\}$ travelers associated with a predefined starting and ending depot. The profit from a set can be collected if precisely one node of a set is visited by a traveler.

\begin{figure}[htbp]
    \centering
    \includegraphics[width=0.5\textwidth]{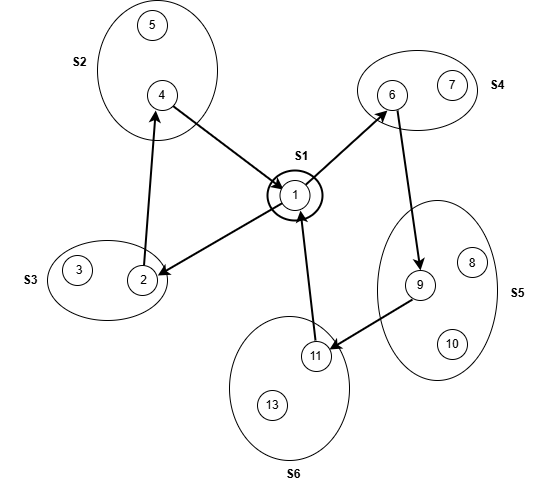} 
    \caption{An example of the sDmSOP.}
    \label{fig:sDmSOP}
\end{figure}

Here, Figure \ref{fig:sDmSOP} illustrates the example of the solution of the sDmSOP instance using two travelers. The set $s_1$ represents the depot, and from there two travelers start and end their journey. The traveler $t_1$ visits set $s_3$ and $s_2$, The traveler $t_2$ visits $s_4$, $s_5$ and $s_6$ respectively. 

\subsection{Mathematical Formulation for the sDmSOP}
\label{subsection:2.1}

To represent an ILP formulation for the sDmSOP, we use some notations, which are given as follows:

\begin{itemize}

  \item  $t \in \{\, 1, \ldots, m \,\} $ (the different salesman).
  
  \item $i, j \in \{\, 1, \ldots, n \,\} $ (the list of vertices).
  
  \item The depot is at vertex $1$.
  
  \item $c_{ij}$ represents the edge weight of the edge $(i, j)$.
  
  \item  $P_q$ represents the profit associated with a set $s_q$.
  
  \item $B$ represents the budget.
  
\end{itemize}

An ILP formulation is constructed using the decision variables:

\begin{itemize}

    \item $x_{tij} = 1$ if the traveler $t$ uses the edge $(i,j) \in E $ and $0$ otherwise.
    
    \item $y_{ti} = 1$ if the vertex $i$ is visited by the traveler $t$ and $0$ otherwise.
    
    \item  $z_{tq} = 1$ if any vertex in the set $s_q$ is visited by the traveler $t$ and $0$ otherwise.
    
    \item  $u_{ij}$ is the flow variable for the Sub-tour Elimination Constraints (SECs).
    
\end{itemize}

The proposed mathematical formulation of the sDmSOP is as follows:

\begin{equation}
\label{eq:1}
\text{maximize} \: \sum_{t} \sum_{q} P_q z_{tq}, 
\end{equation}

subject to:

\begin{equation}
\label{eq:2}
x_{tij},y_{ti},z_{tq} \in \{\, 0,1 \,\}, \quad \forall t, \forall i, \forall j, \forall q,
\end{equation}

\begin{equation}
\label{eq:3}
 \sum_{i} \sum_{j} x_{tij} c_{ij} \leq B, \quad \forall t, 
\end{equation}

\begin{equation}
\label{eq:4}
 \sum_{t} \sum_{j} x_{t1j} = m = \sum_{t} \sum_{j} x_{tj1}, 
\end{equation}

\begin{equation}
\label{eq:5}
 \sum_{i \in V-\{\, j \,\}}  x_{tij}=y_{tj}, \ \quad \forall t, \forall j, 
\end{equation}

\begin{equation}
\label{eq:6}
 \sum_{i \in V-\{\, j \,\}}  x_{tji}=y_{tj}, \ \quad \forall t, \forall j, 
\end{equation}

\begin{equation}
\label{eq:7}
 \sum_{i \in s_q}  y_{ti}=z_{tq}, \ \quad \forall t, \forall q, 
\end{equation}

\begin{equation}
 \sum_{t}  z_{tq}\leq 1, \ \quad \forall q, \label{eq:8}
\end{equation}

\begin{equation}
\label{eq:9}
 0 \leq u_{ij} \leq (n-m) \sum_{t=1}^{m}x_{tij}, \ \quad \forall i, \forall j, 
\end{equation}

\begin{equation}
\label{eq:10}
 \sum_{j \in V}u_{ij} -\sum_{j \in V-\{\, 1 \,\}}u_{ji}=\sum_{t=1}^{m}y_{ti}, \ \quad \forall i \: , \: i \neq j, i \in V-\{\, 1 \,\}. 
\end{equation}

The objective function \eqref{eq:1} maximizes collected profits from the sets visited, constraints \eqref{eq:2} define the domain of the variables $x_{tij}, y_{ti} \text{, and}\, z_{tq}$. Constraint \eqref{eq:3} ensures that each traveler's budget $B$ is not exceeded. Constraint \eqref{eq:4} ensures that exactly $m$ travelers start and end at depot 1. Constraints \eqref{eq:5} and \eqref{eq:6} imply that the in-degree is equal to the out-degree of a vertex except for the depot. Constraint \eqref{eq:7} ensures that a set $s_q$ is visited by a traveler $t$ if any vertex in the set is visited and at most one vertex can be visited per set. Constraint \eqref{eq:8} implies that no set can be visited by more than one traveler while constraints \eqref{eq:9} and \eqref{eq:10}, are used to remove the sub-tours in the path. The SECs are based on the TSP model proposed by \cite{gavish1978travelling} and assessed for the Asymmetric Traveling Salesman Problem (ATSP) by \cite{oncan2009comparative}.

In the above ILP formulation, equations \eqref{eq:1}-\eqref{eq:10} attempt to find out the optimal path with maximization of the profit using permutation of the sets and the vertices which are to be visited in the specific set.

The following two sections discuss implementing the GA and VNS to solve the sDmSOP.

\section{Genetic Algorithm for the sDmSOP}
\label{section:3}

Genetic Algorithm (GA) is a meta-heuristic optimization technique inspired by the principles of natural selection and genetics first proposed by \cite{holland1975adaptation}. It is being used to solve many complex problems in engineering and computer science. Recently, the GA and modified version of GA have been used to solve the TSP and related variants by \cite{george2020genetic}, \cite{prayudani2020analysis}, \cite{ilin2023hybrid}, \cite{zheng2023reinforced}, and \cite{carrabs2021biased}.

The algorithm initiates with an initial population comprising candidate solutions, depicted as chromosomes. These chromosomes then undergo a sequence of genetic operators, including selection, crossover, and mutation, to generate a subsequent population consisting of offspring chromosomes. This process mimics natural evolution, where the fittest individuals are more likely to survive and reproduce, passing favorable traits to the next generation. The fitness of each chromosome is evaluated using an objective function, and the best-fit individuals are selected to form the next generation. The algorithm continues to iterate through the generations until a stopping criterion of reaching a maximum number of generations or achieving a satisfactory fitness level is met.

In the following section, operations used in the GA are explained thoroughly, while Algorithm \ref{Algorithm 1} represents the general scheme depicted by the GA. It starts with initializing the population size $P$ (i.e., $200$). Then, it evaluates the population's fitness (i.e., profit). The loop runs until the stopping condition is met, and our algorithm stops if it cannot find a better solution in $50$ iterations. 

 \begin{algorithm}[h!]
    \SetAlgoLined
    \KwResult{Returns the Valid Path}
    $P \leftarrow$ \textit{InitializePopulation ()}\\
    \textit{EvaluateFitness ($P$)}\\
    \While {termination$\_$condition is false} {
        $\eta \leftarrow $ \textit{Selection ($P$)}\\
        $\eta \leftarrow $ \textit{Crossover ($\eta$)}\\
        $\eta \leftarrow $ \textit{Mutation ($\eta$)}\\
        \textit{EvaluateFitness ($\eta$)}\\
        $P \leftarrow$ \textit{ReplacePopulation ($P$, $\eta$)}\\
    }
    \caption{Genetic Algorithm}
    \label{Algorithm 1}
\end{algorithm}

\subsection{Chromosome Representation}
\label{subsection:3.1}

\begin{itemize}

  \item  Each chromosome has the same length.
  
  \item  Each chromosome stores the information about clusters traversed by each salesman, where the number of salesmen is specified during input.
  
  \item A chromosome also stores information about whether a cluster should be part of the final output. This means that every cluster is not selected in the problem, and the final answer can consist of a subset of clusters instead of all. 
  
\end{itemize}

\begin{figure}[htbp]
    \centering
    \includegraphics[width=0.35\textwidth, height=0.15\textwidth]{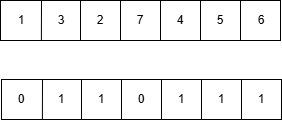}
    \caption{Arrangement and Membership array in Chromosome.}
    \label{fig: Chromosome Representation}
\end{figure}

Fig. \ref{fig: Chromosome Representation} illustrates the representation of a chromosome in the example explained below. A chromosome consists of two arrays: Arrangement and Membership. Both arrays have the same size of  ($\text{cluster count} + \text{salesman count} - 1$). The Arrangement Array (AA) initially stores the indices of the clusters in random order. The elements in an AA are of two categories: A cluster number and a separator. The cluster number represents the cluster we target, and the separator indicates that a new salesman is entering the picture. The number is treated as a separator if its value exceeds the cluster count. The cluster $1$ is ignored as it is the starting and ending cluster for all the travelers.

For more clarity, let us take an example based on Fig. \ref{fig:sDmSOP}, where the cluster count is $6$, and the number of salesmen used is $2$. Here, the array length for the chromosome representation will be $7$ (i.e. $6+2-1$). Consider a random AA: $[1, 3, 2, 7, 4, 5, 6]$. The number used as separator is $7$, as it exceeds the cluster count. Hence, the salesmen will be assigned cities $[[1, 3, 2], [4, 5, 6]]$. Being a depot, the cluster $1$ is ignored; this means the array becomes $[[3, 2], [4, 5, 6]]$  where the salesman $1$ is given the cluster $3$ and $2$, and the salesman $2$ is given the clusters $4$ and $5$, and $6$ respectively.
 
The Membership Array (MA) keeps track of which cluster to include in fitness calculation. It is a boolean array. Parameter one\_rate specifies the probability of a given cluster being part of the final fitness calculation during initialization. Membership of separators and city one is ignored. 

The combination of the membership array $[0, 1, 1, 0, 1, 1, 1]$ and the arrangement array $[1, 3, 2, 7, 4, 5, 6]$ results in the cluster sequence $[1\to3\to2\to1, 1\to4\to5\to6\to1]$.

\subsection{The \textit{Selection} Operation}
\label{subsection:3.2}

In a genetic algorithm, the \textit{Selection} operation is responsible for choosing individuals from the previous generation to ensure the population's quality. 

The probability of being chosen is directly proportional to the chromosome's fitness value, which means that chromosomes with a better fitness score have a higher chance of being selected. Here, selection is based only on the chromosomes' fitness value (total profit). By selecting the fittest individuals, the genetic algorithm can gradually improve the population's quality and, over time, converge towards a globally optimum solution. 
  
\subsection{Fitness Calculation for a Chromosome}
\label{subsection:3.3}

\begin{itemize}

    \item The cluster sequence for each salesman is extracted from a given chromosome.
    
    \item Utilizing dynamic programming, the optimal cities to be visited for each cluster sequence are calculated. This process determines the path from the initial cluster ($1$) to the final cluster ($1$) with the minimum distance. Notably, this information is recalculated each time the function encounters a new cluster sequence and is not permanently stored.
    
    \item Fitness is defined as the sum of profits from each cluster sequence.
    
    \item If the salesman count does not match the specified count or the cost constraint is unspecified for any chromosome, the fitness is returned as $0$.
    
\end{itemize}

\subsection{The \textit{Crossover} Operation}
\label{subsection:3.4}

\begin{figure}[htbp]
    \centering
    \includegraphics[width=0.40\textwidth, height=0.25\textwidth]{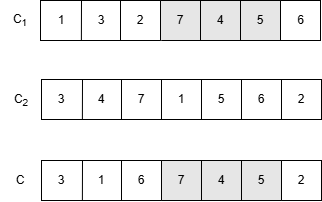} 
    \caption{Crossover.}
    \label{fig:Crossover}
\end{figure}

The \textit{Crossover} operation is a crucial genetic algorithmic procedure used to generate a new offspring chromosome from two parent chromosomes, $c_1$ and $c_2$, by exchanging genes between selecting a random region within the chromosomes and exchanging genes within the particular region creates a new offspring chromosome. Fig. \ref{fig:Crossover} represents how the new offspring is generated by the GA using crossover. It first selects a random region within the chromosomes and identifies the genes within the region of the first chromosome $c_1$. The \textit{Crossover} operation is performed by copying the genes in $c_2$ not in the crossover region of $c_1$ to the beginning of a new chromosome $c$. The crossover region of $c_1$ is then copied to $c$, and the remaining genes of $c_2$ that were not already in $c$ are copied to the end of $c$. The resulting offspring chromosome has the same length as the parent chromosomes and contains unique genes. The \textit{Crossover} operation thus effectively generates new offspring for the genetic algorithm.

\subsection{The \textit{Mutation} Operation}
\label{subsection:3.5}

In this operation, the function takes a mutation rate as input, a 5\% probability value that determines the likelihood of a gene in the chromosome undergoing mutation. The function first iterates over each gene in the arrangement vector of the chromosome and randomly swaps it with another gene with a probability determined by the mutation rate. This step shuffles the order of the genes in the chromosome to generate new permutations. The \textit{frand} function returns a random float value between 0 and 1, which is compared to the mutation rate to determine if the gene is selected for mutation. The function then iterates over each element in the membership vector of the chromosome and flips its Boolean value with a probability determined by the mutation rate. This step toggles the membership value of the corresponding gene to its opposite value. The \textit{frand} function is used again to determine if the element is selected for mutation.

In summary, the \textit{Mutation} operation randomly mutates the genetic material of a chromosome by shuffling the order of its genes and flipping the membership value of some of its elements. This introduces randomness and diversity in the population and helps the genetic algorithm avoid being stuck in local optima.\\

\subsection{Procedure to Find the Cities with the Minimum Distance in the Given Sequence of Clusters}
\label{subsection:3.6}

Employing 2D dynamic programming, consider the cluster sequence: $1\to3\to2\to1$ based on Fig. \ref{fig:sDmSOP}. In this sequence, cluster 3 includes cities 2 and 3, while cluster 2 comprises cities 4, and 5.

\begin{itemize}

    \item Visualize this as a neural network with four layers. Layer $1$ has the city $1$. Layer $2$ has the cities $2$ and $3$. Layer $3$ has the cities $4$, and $5$. Layer $4$ has the city $1$. Each layer is fully connected to the next layer.
    
    \item $dp[i,j]$ stores the minimum distance to reach the city present at the index $j$ of the layer $i$. In the above example, $dp[1,]$ is an array of size $1$, $dp[2,]$ is an array of size $2$, $dp[3,]$ is an array of size $2$, and $dp[4,]$ is an array of size $1$. $dp[i,]$ stores the minimum costs to reach each city in layer $i$, and each element in this array is derived by considering all possible paths from the previous layer. 
    
    \item $dp[3,1]$ means the minimum distance to reach the city $4$ ($dp[3,]$ corresponds to the cities $4$ and $5$, at index $1$, we have the city $4$, and at index $2$, the city $5$), that is present in cluster indexed at $3$ (which is the cluster $2$) from the first layer, which consists of the only city $1$. 
    
\end{itemize}

We then apply the following formula to calculate $dp[i,j]$:

\begin{itemize}

    \item $dp[1,1] = 0$,
    
    \item $dp[i, j] = \min_{k} \left( dp[i-1, k] + \text{dis}(\text{city}_{i-1}[k], \text{city}_{i}[j]) \right)$.
    
\end{itemize}

where:

\begin{itemize}

\item \( dp[i-1, k] \) is the minimum cost to reach the city at index \( k \) in layer \( i-1 \),

\item \(\text{dis}(\text{city}_{i-1}[k], \text{city}_{i}[j])\) is the distance or cost to travel from the city at index \( k \) in layer \( i-1 \) to the city at index \( j \) in layer \( i \).

This equation ensures that \( dp[i, j] \) holds the minimum cost required to reach the city at index \( j \) in the $i$th layer from any city in the $(i-1)$th layer.

\end{itemize}

\section{Variable Neighborhood Search for the sDmSOP}
\label{section:4}

Variable Neighborhood Search (VNS) is a robust metaheuristic extensively utilized for tackling combinatorial optimization problems, notably the Set Orienteering Problem (SOP). Introduced by \cite{mladenovic1997variable}, VNS iteratively refines solutions through diverse neighborhood operators, thereby enhancing the exploration of solution spaces. The approach initiates with a greedy solution and aims to minimize the distance per additional profit from visiting new clusters. Subsequently, the VNS algorithm systematically applies predefined neighborhood operators to iteratively improve the current best solution and obtain the best approximation for the problem.

The application of VNS to the Orienteering Problem (OP) was further developed by \cite{sevkli2006variable}, who introduced the first algorithm solely based on VNS to solve the OP. Their work demonstrated the effectiveness of VNS in finding high-quality solutions to the OP, showcasing its potential to escape local minima and explore the solution space more thoroughly. Subsequent research, such as that by \cite{pvenivcka2017dubins(a)}, \cite{pvenivcka2017dubins(b)}, and \cite{pvenivcka2019variable} has further explored the use of VNS and similar neighborhood operators for initial solutions in other variants of the OP, including the Dubins Orienteering Problem (DOP) and the Orienteering Problem with Neighborhoods (OPN).

VNS employs a dual strategy of \textit{local\_search} and \textit{shake} to refine an initial solution systematically. The iterative process aims to optimize path profits within a local neighborhood, pushing towards a superior local maximum with minimal cost. It is the essence of the \textit{local\_search} procedure. Conversely, the \textit{shake} operation introduces controlled randomness by adding or substituting clusters, broadening the exploration to escape potential premature local maxima.

In practice, VNS continually alternates between these procedures, generating and refining solutions immediately to improve performance. During each iteration, if a new solution is valid and surpasses the current one in profit while adhering to the constraints, the algorithm updates and signals an improvement. This careful balance of local refinement and randomized exploration distinguishes VNS, ensuring a dynamic search that efficiently navigates solution space.

\begin{algorithm}[h!]
    \SetAlgoLined
    \KwResult{$u$ - Returns the Optimal Solution}
    $u \leftarrow$ \textit{constructInitialSolution()}\\ 
    \While {true} {
        $l \leftarrow 1$\\
        \While {$l \leq l_{max}$} {
            $u' \leftarrow$ \textit{shake(u, l)}\\
            $u'' \leftarrow$ \textit{local\_search($u', l$)}\\
            \eIf{\textit{isValid($u''$)} \textbf{and} $P(u'') > P(u)$} {
                $u \leftarrow u''$\\
                $l \leftarrow 1$\\
            } {
                $l \leftarrow l + 1$\\
            }
        }
    }
    \caption{Variable Neighborhood Search}
\end{algorithm}

\subsection{The \textit{constructInitialSolution} Function}
\label{subsection:4.1}

Generating a computationally light and viable initial solution is a crucial first step for the problem, and it is more challenging due to the multiple travelers associated with the single depot. 

So, we used a construction procedure of the initial solution using a greedy approach that minimizes the additional length of the path per additional profit using the following formula: We iterate over all the sets, and for each set, we find a position for it in the sequence for each salesman. Essentially, we iterate over all the positions and choose the position $\text{pos}$ such that: 
  \[\text{pos} = \arg \min_{j} \frac{L(\text{adj}_{ij})-L(\text{adj})}{p_i},\] 
  where $p_i$ is the profit from the newly added set, without violating the maximum cost ($B$) condition,  $L(\text{adj}_{ij})$ is the cost of the configuration after adding set $i$ at position $j$ (from the start) and $L(\text{adj})$ is the cost associated to the old arrangement. Then, a Hungarian Assignment Algorithm is used to improve the minimal-cost assignment for each traveler. 
 
This method provided a straightforward mechanism to check the validity of the case, relying on two conditions: 

\begin{enumerate}

\item The total number of travelers $m$ should be, at most, the number of sets $p$.

\item The distance bound for each traveler should be less than the assigned budget of $B$.

\end{enumerate}

This solution is then optimized to construct a valid solution using VNS.

It iterates over each unvisited cluster and explores potential positions in the paths to insert the unvisited cluster. At each iteration, it evaluates the validity and optimizes the solution if a better position is found. Sections \ref{subsection:4.2} and \ref{subsection:4.3} explain the essential optimizers used to solve the sDmSOP.

\subsection{The \textit{shake} Operation}
\label{subsection:4.2}

The \textit{shake} operation employs two procedures: \textit{Path Move} and \textit{Path Exchange} to introduce randomized changes to the solution vector for efficient exploration of new neighborhoods. The process dynamically modifies the solution vector to adjust the traversal order of clusters and potentially introduce new, unvisited clusters. Strict adherence to hard constraints is maintained, retaining only new values that satisfy all the constraints.\\

\noindent\textbf{\textit{Path Move}}

The \textit{Path Move} procedure introduces randomness by selecting a path, represented as a subarray and an index. The path chosen is then relocated to the new index, and the cost is evaluated. It is crucial to note that the same salesman must take the entire path.

\begin{enumerate}

    \item[{\bfseries 1.}]\textbf{Selection of Neighborhoods and Positions}
    
    \begin{itemize}
    
         \item Random non-empty neighborhoods, \( r_{\text{Swap}} \) and \( r_{\text{Range}} \), are selected from the existing solution. 
         
         \item \( r_{\text{Swap}} \) represents the traveler to which the new path will be moved, and \( r_{\text{Range}} \) represents the traveler from which the path to be moved.
         
        \item A random position, \(id_{\text{Swap}} \), within \( r_{\text{Swap}} \) is chosen to determine the index around which the swapped clusters will be inserted. A random choice (\textit{whichSide}) determines whether the swapped clusters will be placed before or after \( id_{\text{Swap}} \). 
        
        \item A range \( [lo, hi] \) within \( r_{\text{Range}} \) is selected, where \(lo\) and \(hi\) represent indexes in the path of \( r_{\text{Range}} \) where clusters will be swapped out.
        
    \end{itemize}

    \item[{\bfseries 2.}]\textbf{Performing the Swap}
    
    \begin{itemize}
    
        \item If \( r_{\text{Swap}} \) equals \( r_{\text{Range}} \), that is, the traveler to and from which the part is to be moved are the same, the algorithm ensures that the swap operation does not violate the selected range \( [lo, hi] \) within the same neighborhood.
        
        \item If \( r_{\text{Swap}} \) differs from \( r_{\text{Range}} \), that is, if the traveler to and from the part is to be moved are different, the algorithm inserts elements in the range of indexes \( [lo, hi] \) from \( r_{\text{Range}} \) to \textit{whichSide} (before or after) of \(id_{\text{Swap}} \) for traveler \( r_{\text{Swap}} \).
        
        \item The modified neighborhoods are stored in \( u' \).
        
    \end{itemize}
    
\end{enumerate}

\noindent\textbf{\textit{Path Exchange}}

The \textit{Path Exchange} procedure introduces randomness by selecting two paths, each represented as a subarray, originating from distinct salesmen, and exchanging them. This operation involves relocating entire paths between different salesmen.

\begin{enumerate}

    \item[{\bfseries 1.}] \textbf{Selection of Neighborhoods and Positions}
    
        \begin{itemize}
        
            \item Two random non-empty neighborhoods, \( r_{\text{P1}} \) and \( r_{\text{P2}} \), are selected from the solution.
            
            \item Random positions within each neighborhood, represented by \({id}[0] \), \({id}[1] \) for \( r_{\text{P1}} \) and \({id}[2] \), \({id}[3] \) for \( r_{\text{P2}} \), are chosen to determine the range of clusters that will be swapped.
            
        \end{itemize}
        
    \item[{\bfseries 2.}] \textbf{Performing the Swap}
    
        \begin{itemize}
        
            \item If \( r_{\text{P1}} \) and \( r_{\text{P2}} \) are the same neighborhood, additional checks are performed to ensure that the selected positions are different.
            
            \item Clusters within the specified ranges are swapped between the two neighborhoods.
            
            \item The modified neighborhoods are stored in \( u'\).
            
        \end{itemize}
        
\end{enumerate}

\subsection{The \textit{local\_search} Operation}
\label{subsection:4.3}

The \textit{local\_search} operation within the VNS algorithm employs two neighborhood structures to optimize the current solution vector. The \textit{local\_search} operation consists of two procedures, \textit{One Cluster Move} and \textit{One Cluster Exchange} to increase randomness in the solution.\\

\noindent\textbf{The \textit{One Cluster Move} Procedure}

The \textit{One Cluster Move} procedure aims to refine the current solution vector \(u\) by selecting two random positions, \(i_1\) and \(i_2\), and exploring two possible arrangements: placing \(i_1\) after \(i_2\) or placing \(i_2\) before \(i_1\). The modified solution is retained only if it increases the profit (\(P\)).

\begin{enumerate}

    \item[{\bfseries 1.}]\textbf{Initialization and Cluster Selection}
    
        \begin{itemize}
        
            \item Generate random indexes, $i$ and $j$, from the range of visited paths.
            
            \item If either $i$ or $j$ exceeds the total number of available paths, or if they are equal, revert to the original configuration.
            
            \item Locate the positions of the clusters represented by $i$ and $j$ within the solution configuration.
            
        \end{itemize}
        
    \item[{\bfseries 2.}]\textbf{Random Choice and Path Modification}
    
        \begin{itemize}
        
            \item Generate a random binary value of $1$ or $0$.
            
            \item If the value is $1$, relocate the cluster $i$ after the cluster $j$ in their respective paths, removing $i$ from its original position.
            
            \item If the value is $0$, relocate the cluster $j$ before the cluster $i$ in their respective paths, removing $j$ from its original position.
            
        \end{itemize}

\end{enumerate}

\noindent\textbf{The \textit{One Cluster Exchange} Procedure}

The \textit{One Cluster Exchange} procedure refines the current solution vector \(u\) by selecting two random clusters, \(i\) and \(j\), and exchanging them. The modified solution is obtained by swapping the positions of clusters \(i\) and \(j\) within the solution vector \(u\).

\begin{enumerate}

    \item[{\bfseries 1.}]\textbf{Initialization and Cluster Selection}
    
        \begin{itemize}
        
            \item Generate two random indexes, $i$ and $j$, within the range of the visited clusters.
            
            \item Ensure that neither $i$ nor $j$ exceeds the total number of available paths.
            
        \end{itemize}
        
    \item[{\bfseries 2.}]\textbf{Locate and Swap Clusters}
    
        \begin{itemize}
        
            \item Locate the positions of the clusters represented by $i$ and $j$ within the solution configuration using the indexes.
            
            \item Swap the clusters $i$ and $j$ in their respective paths and check the validity of the solution.
            
        \end{itemize}
        
\end{enumerate}

\section{Computational Tests}
\label{section:5}

The mathematical model is simulated on GTSP instances using GAMS 37.1.0 with CPLEX, and the performance comparison with multiple travelers is presented. The simulation is performed on a Windows 10 platform with an i7-6400 CPU @3.4Ghz processor and 32GB of RAM. The generation of sDmSOP instances is detailed in section \ref{subsection:5.1}, and the simulation results are provided in section \ref{subsection:5.2}.

\subsection{The Test Instances}
\label{subsection:5.1}

The Generalized Traveling Salesman Problem (GTSP) instances of \cite{noon1988generalized} are used to evaluate the comparative results of the mathematical formulation and algorithms.  

To adapt the GTSP instances for the sDmSOP, the following modifications are done:

\begin{itemize}

    \item Remove node $1$ from the original cluster and move it to a new cluster (i.e., depot) that contains only node $1$. This new cluster is the starting and ending point for every traveler.

The profit is generated for each of the clusters based on the following rules:

\item Profit function $g_1$:

    \begin{enumerate}
    
        \item Cluster profit $=$ number of nodes in the cluster.
        
        \item Profit of cluster $1$ is $0$.
        
    \end{enumerate}
    
\item    Profit function $g_2$:

    \begin{enumerate}
    
        \item Cluster profit $=$ Sum of profits of nodes of each cluster.
        
        \item The $i$th node profit $ = (1+(7141\times{i}))\: \text{ mod } 100$.
        
    \end{enumerate}

\end{itemize}

\subsection{Computational Results}
\label{subsection:5.2}

In this section, we present the simulation results of the ILP, GA, and VNS when $w=0.25$ in Table \ref{tab: Table 1} as follows: The first four columns represent the GTSP instance used, number of nodes ($n$) in the instance, number of travelers ($t$) used to solve the instance, the profit rule ($P_g$) to calculate the profit of a cluster, the following six columns represent the solution and time (in seconds), given by the ILP, GA and VNS respectively. 

Table \ref{tab: Table 2} is arranged as follows: The first four columns represent the Set number (i.e., Set $1$, Set $2$ or Set $3$), Instance name, $n$ represents the number of nodes available in the instance, $t$ represents the number of travelers used. The following eight columns represent the solution and time (in seconds) provided by the VNS and GA for the profit rules $g_1$ and $g_2$. Set $1$ contains the instances that have less than $200$ nodes, Set $2$ contains the instances in which the number of nodes is in the range of $200-500$, and Set $3$ has all the instances that have more than $500$ nodes to $217vm1084$, which is the largest instance available in the GTSP instance. 

In Table \ref{tab: Table 1}, small GTSP instances of less than 100 nodes are taken for simulation. CPLEX is not able to solve $5$ instances optimally and goes Out Of Memory (OOM); the rest of the $15$ instances are solved optimally. GA and VNS take very little time to solve these instances compared to CPLEX. We observe that if the number of travelers is two, GA takes less time than VNS, but as number of travelers increases from two to three, GA takes too much time to produce the same results as VNS. 

\begin{longtable}[h!]{cccccccccccc}
\caption{Comparison of ILP, VNS, and GA on the GTSP data instances of the sDmSOP when ($w = 0.25$)}
\label{tab: Table 1} \\
\hline \endhead
Instance   & $n$  & $t$ & $P_g$ & \multicolumn{2}{c}{ILP}     &  & \multicolumn{2}{c}{GA}  &  & \multicolumn{2}{c}{VNS}          \\ 
\cline{5-6}\cline{8-9} \cline{11-12}
           &    &            &    & Opt. Solution & Time (sec.) &  & Solution & Time (sec.) &  & Solution & Time (sec.) \\ 
\hline

11berlin52	 &52	&2	&$g_1$	&37	          &114.69	    &&37	&1.82	    &&37	   
&3.92\\
11eil51	     &51	&2	&$g_1$	    &24	      &293.22	    &&24	&2.45	    &&24	   &2.95\\
14st70	     &70	&2	&$g_1$	    &OOM	  &12777.98     &&27	&8.42	    &&27	   &3.47\\
16eil76	     &76	&2	&$g_1$	    &OOM	  &525118.19    &&40	&4.35	    &&40	   &4.21\\
11berlin52	 &52	&3	&$g_1$	    &37	      &245.73	    &&37	&274.29 	&&37	   &4.39\\
11eil51	     &51	&3	&$g_1$	    &28	      &992.53	    &&28	&153.87 	&&28	   &3.78\\
14st70	     &70	&3	&$g_1$	    &27	      &1848155.81	&&27	&455.15 	&&27	   &3.99\\
16eil76	     &76	&3	&$g_1$	    &OOM	  &85308.55 	&&45	&164.52 	&&45	   &5.19\\
11berlin52	 &52	&2	&$g_2$	    &1729	  &240.17	    &&1729	&2.11	    &&1729	&4.12\\
11eil51	     &51	&2	&$g_2$	    &1279	  &250.41	    &&1279	&1.39	    &&1279	&2.87\\
14st70	     &70	&2	&$g_2$	    &OOM	  &26188.91 	&&1271	&9.08	    &&1271	&3.40\\
16eil76	     &76	&2	&$g_2$	    &2192	  &594952.39	&&2192	&6.48	    &&2192	&4.21\\
11berlin52   &52	&3	&$g_2$	    &1729	  &261.08	    &&1729	&286.34 	&&1729	&4.33\\
11eil51	     &51	&3	&$g_2$	    &1466	  &819.53	    &&1466	&154.43 	&&1466	&3.82\\
14st70	     &70	&3	&$g_2$	    &1271	  &2019069.03	&&1271	&465.10 	&&1271	&4.09\\
16eil76	     &76	&3	&$g_2$	    &OOM	  &1319603.16	&&2394	&218.18 	&&2394	&5.15\\

\hline
\end{longtable}

\begin{longtable}[h!]{cccclcclcccccccccc}
\caption{Profit comparison of VNS and GA on the GTSP data instances of the sDmSOP when ($w=0.25$)}
\label{tab: Table 2} \\
\hline \endhead
              &            &      &            &  & \multicolumn{2}{c}{$g_1$}               &  & \multicolumn{2}{c}{$g_2$}  &  & \multicolumn{2}{c}{$g_1$}               &  & \multicolumn{2}{c}{$g_2$}                 \\ 
\cline{6-7}\cline{9-10} \cline{12-13}\cline{15-16}
              & Instance   & $n$    & $t$ &  & \multicolumn{2}{c}{VNS}              &  & \multicolumn{2}{c}{VNS}  &  & \multicolumn{2}{c}{GA}              &  & \multicolumn{2}{c}{GA}                \\
              &            &      &            &  & Solution          & Time                  &  & Solution           & Time        &  & Solution          & Time             &  & Solution           & Time         \\ 
\hline
Set 1 
&11berlin52	&52	    &2	&&37	&3.92	&&1729	&4.12	&&37	&1.92	&&1729	&1.83\\
&11eil51	&51	    &2	&&24	&2.95	&&1279	&2.87	&&24	&2.18	&&1279	&1.84\\
&14st70	    &70	    &2	&&27	&3.47	&&1271	&3.40	&&27	&8.30	&&1271	&8.18\\
&16eil76	&76	    &2	&&40	&4.21	&&2192	&4.21	&&40	&5.02	&&2192	&6.65\\
&20kroA100	&100	&2	&&39	&4.43	&&1925	&4.49	&&37	&5.05	&&1755	&5.14\\
&20kroB100	&100	&2	&&45	&4.64	&&2251	&4.59	&&44	&4.66	&&2251	&5.53\\
&20kroC100	&100	&2	&&34	&4.25	&&1597	&4.14	&&34	&8.21	&&1597	&6.39\\
&20kroD100	&100	&2	&&37	&4.34	&&1594	&4.28	&&37	&7.73	&&1594	&6.99\\
&20kroE100	&100	&2	&&49	&4.88	&&2388	&4.86	&&45	&5.22	&&2267	&6.48\\
&20rat99	&99	    &2	&&30	&3.91	&&1470	&3.88	&&30	&17.43	&&1470	&13.62\\
&20rd100	&100	&2	&&38	&4.50	&&1886	&4.41	&&38	&11.38   &&1797	&7.31\\
&21eil101	&101	&2	&&57	&4.97	&&2969	&5.16	&&57	&3.75	&&2969	&3.90\\
&21lin105	&105	&2	&&29	&4.56	&&1436	&4.60	&&29	&21.94	&&1436	&19.15\\
&22pr107	&107	&2	&&29	&4.66	&&1453	&4.66	&&29	&15.04	&&1453	&15.08\\
&25pr124	&124	&2	&&41	&4.87	&&1977	&4.85	&&37	&9.67	&&1785	&16.09\\
&26bier127	&127	&2	&&108	&8.40	&&5196	&8.49	&&102	&3.61	&&4921	&3.37\\
&26ch130	&130	&2	&&57	&5.54	&&3215	&5.64	&&55	&4.51	&&3215	&4.01\\
&28pr136	&136	&2	&&38	&5.31	&&1817	&5.30	&&38	&12.99	&&1817	&19.75\\
&29pr144	&144	&2	&&54	&6.19	&&2724	&6.29	&&48	&9.45	&&2410	&9.85\\
&30ch150	&150	&2	&&49	&5.24	&&2361	&5.17	&&43	&8.13	&&2361	&13.08\\
&30kroA150	&150	&2	&&49	&5.49	&&2533	&5.51	&&47	&6.15	&&2444	&10.24\\
&30kroB150	&150	&2	&&60	&6.18	&&2732	&6.22	&&53	&9.70	&&2634	&6.70\\
&31pr152	&152	&2	&&37	&4.90	&&1804	&4.89	&&33	&23.64	&&1699	&17.96\\
&32u159	    &159	&2	&&55	&6.51	&&3303	&6.77	&&54	&13.57	&&3228	&14.75\\
&39rat195	&195	&2	&&53	&5.84	&&2634	&5.78	&&44	&19.32	&&2425	&40.38\\
            
\textbf{Avg.} &            &      &            &  & \textbf{44.64}   & \textbf{4.97}  &  & \textbf{2229.44}  & \textbf{4.98}       &  & \textbf{42.48}   & \textbf{9.54}      &  & \textbf{2159.96	}   & \textbf{10.57} \\

&11berlin52 &52 &3 &&37 &4.38 &&1729 &4.33 &&37 &248.75 &&1729 &254.53\\
&11eil51 &51 &3 &&28 &3.78 &&1466 &3.82 &&28 &145.59 &&1466 &145.27\\
&14st70 &70 &3 &&27 &3.99 &&1271 &4.09 &&27 &425.10 &&1271 &423.41\\
&16eil76 &76 &3 &&45 &5.19 &&2394 &5.14 &&45 &139.25 &&2346 &106.19\\
&20kroA100 &100 &3 &&51 &5.66 &&2554 &5.65 &&48 &79.28 &&2131 &77.16\\
&20kroB100 &100 &3 &&50 &5.40 &&2489 &5.32 &&50 &52.96 &&2489 &70.30\\
&20kroC100 &100 &3 &&39 &4.82 &&1974 &4.89 &&39 &154.23 &&1974 &194.34\\
&20kroD100 &100 &3 &&39 &4.71 &&1640 &4.70 &&37 &188.18 &&1594 &151.95\\
&20kroE100 &100 &3 &&55 &5.56 &&2554 &5.59 &&52 &160.71 &&2249 &113.38\\
&20rat99 &99 &3 &&30 &4.47 &&1470 &4.57 &&30 &458.14 &&1470 &545.77\\
&20rd100 &100 &3 &&38 &5.02 &&1886 &4.95 &&38 &242.85 &&1886 &198.37\\
&21eil101 &101 &3 &&69 &6.25 &&3415 &6.10 &&66 &32.12 &&3375 &35.64\\
&21lin105 &105 &3 &&29 &5.16 &&1436 &5.17 &&29 &984.53 &&1436 &1011.38\\
&22pr107 &107 &3 &&29 &5.26 &&1453 &5.22 &&29 &542.73 &&1453 &544.62\\
&25pr124 &124 &3 &&41 &5.41 &&1977 &5.35 &&41 &232.38 &&1977 &220.56\\
&26bier127 &127 &3 &&112 &8.82 &&5495 &8.75 &&108 &20.57 &&5142 &13.78\\
&26ch130 &130 &3 &&62 &6.31 &&3302 &6.38 &&63 &22.86 &&3735 &48.99\\
&28pr136 &136 &3 &&41 &5.28 &&2003 &5.30 &&41 &480.39 &&1970 &275.77\\
&29pr144 &144 &3 &&64 &6.65 &&3189 &6.66 &&54 &121.68 &&3189 &177.69\\
&30ch150 &150 &3 &&60 &6.29 &&2887 &6.07 &&56 &124.06 &&2750 &192.00\\
&30kroA150 &150 &3 &&70 &6.67 &&3558 &6.69 &&58 &30.78 &&3017 &33.77\\
&30kroB150 &150 &3 &&74 &7.12 &&3356 &7.02 &&75 &108.48 &&3110 &131.59\\
&31pr152 &152 &3 &&38 &5.24 &&1851 &5.13 &&38 &393.15 &&1851 &481.30\\
&32u159 &159 &3 &&70 &6.70 &&3329 &6.54 &&67 &192.12 &&3568 &186.69\\
&39rat195 &195 &3 &&60 &6.45 &&3008 &6.53 &&48 &307.18 &&2458 &368.16\\

\textbf{Avg.} &            &      &            &  & \textbf{50.32}   & \textbf{5.62}  &  & \textbf{2467.44}  & \textbf{5.60}       &  & \textbf{48.16}   & \textbf{235.52}      &  & \textbf{2385.44}   & \textbf{240.10} \\   

Set 2
&40kroa200  &200 &2 &&77  &6.74  &&3408  &6.66  &&77  &9.50 &&3408 &12.69\\
&40krob200  &200 &2 &&78  &7.83  &&3875  &7.74  &&69  &8.29 &&3081 &6.77\\
&45ts225    &225 &2 &&85  &8.27  &&4228  &8.48  &&55  &14.23 &&3018 &16.96\\
&45tsp225   &225 &2 &&51  &6.42  &&2587  &6.33  &&39  &26.15 &&2105 &52.55\\
&46pr226    &226 &2 &&71  &7.74  &&3634  &7.72  &&55  &14.76 &&2448 &13.43\\
&53gil262   &262 &2 &&67  &7.71  &&3133  &8.50  &&48  &37.68 &&2219 &51.45\\
&53pr264    &264 &2 &&128 &11.46 &&6372  &12.21 &&85  &16.26 &&3413 &11.36\\
&56a280     &280 &2 &&87  &9.46  &&4447  &10.03 &&61  &30.48 &&2691 &26.00\\
&60pr299    &299 &2 &&85  &8.96  &&4130  &9.36  &&52  &47.22 &&2514 &51.14\\
&64lin318   &318 &2 &&132 &12.97 &&6646  &13.02 &&66  &15.32 &&3519 &20.81\\
&80rd400    &400 &2 &&136 &13.14 &&7434  &13.06 &&91  &9.14  &&3853 &8.45\\
&84fl417    &417 &2 &&202 &19.47 &&9954  &19.75 &&94  &33.20 &&4686 &46.24\\
&88pr439    &439 &2 &&221 &18.64 &&10623 &18.73 &&105 &20.20 &&4539 &15.55\\
&89pcb442   &442 &2 &&149 &15.04 &&7382  &14.87 &&65  &85.32 &&3420 &57.79\\

\textbf{Avg.} &&&&& \textbf{112.07} &\textbf{10.99} &&\textbf{5560.93} & \textbf{11.17}        &  & \textbf{68.71}   & \textbf{26.27}      &  & \textbf{3208.14}   & \textbf{27.94} \\

&40kroa200	    &200	&3	&&96	&8.61	&&4835	&8.88	&&78	&63.18	&&3173	&33.46\\
&40krob200	    &200	&3	&&96	&8.98	&&4873	&8.85	&&77	&43.80	&&4087	&37.67\\
&45ts225	    &225	&3	&&93	&9.34	&&4644	&9.03	&&60	&94.05	&&3026	&132.14\\
&45tsp225	    &225	&3	&&59	&7.35	&&2751	&7.20	&&48	&314.77	&&2379	&352.02\\
&46pr226	    &226	&3	&&72	&7.93	&&3696	&7.88	&&64	&119.47	&&3171	&122.58\\
&53gil262	    &262	&3	&&76	&8.93	&&3704	&8.58	&&61	&260.42	&&2462	&420.26\\
&53pr264	    &264	&3	&&131	&12.69	&&6488	&13.10	&&77	&50.08	&&4065	&75.86\\
&56a280	        &280	&3	&&108	&10.60	&&5515	&10.52	&&61	&126.74	&&3156	&182.06\\
&60pr299	    &299	&3	&&97	&10.69	&&4789	&10.15	&&61	&282.72	&&2601	&195.86\\
&64lin318	    &318	&3	&&146	&13.53	&&7423	&12.97	&&76	&68.91	&&4324	&107.08\\
&80rd400	    &400	&3	&&231	&19.90	&&11708	&19.79	&&91	&9.88	&&5250	&16.07\\
&84fl417	    &417	&3	&&217	&21.22	&&10493	&20.74	&&94	&42.21	&&4686	&55.24\\
&88pr439	    &439	&3	&&255	&22.26	&&12506	&22.57	&&105	&40.61	&&4835	&35.21\\
&89pcb442	    &442	&3	&&195	&17.95	&&9706	&17.90	&&72	&93.38	&&3863	&132.95\\

\textbf{Avg.} &&&&& \textbf{133.71} &\textbf{12.86} &&\textbf{6652.21} & \textbf{12.73}        &  & \textbf{73.21}   & \textbf{115.02}      &  & \textbf{3648.43}   & \textbf{135.60} \\

Set 3
&115rat575	&575	&2	&&194	&19.30	&&9508	 &18.51	&&60	&264.29	    &&3279	&227.01\\
&115u574	&574	&2	&&138	&16.22	&&7754	 &16.63	&&66	&149.41	    &&3677	&65.25\\
&131p654	&654	&2	&&190	&23.89	&&9187	 &24.51	&&92	&1486.30	&&4661	&2442.61\\
&132d657	&657	&2	&&150	&16.94	&&6834	 &16.47	&&74	&5070.67	&&3303	&10058.20\\
&145u724	&724	&2	&&219	&22.07	&&10830	 &21.37	&&78	&340.59	    &&3510	&567.21\\
&157rat783	&783	&2	&&229	&23.11	&&11909	 &23.29	&&61	&633.32    	&&3254	&511.57\\
&201pr1002	&1002	&2	&&259	&25.59	&&12161	 &25.98	&&87	&3203.74	&&4463	&3365.56\\
&212u1060	&1060	&2	&&302	&30.24	&&14935	 &30.00	&&91	&2121.65	&&5063	&4919.42\\
&217vm1084	&1084	&2	&&406	&41.76	&&20138	 &41.61	&&122	&860.20	    &&7000	&1158.40\\

\textbf{Avg.} &&&&& \textbf{231.89} &\textbf{24.35} &&\textbf{11472.89} & \textbf{24.26}        &  & \textbf{81.22}   & \textbf{1570.02}      &  & \textbf{4245.56}   & \textbf{2590.58} \\

&115rat575	&575	&3	&&266	&23.79	&&13014	&23.85	&&68	&465.45	    &&3842	&452.70\\
&115u574	&574	&3	&&223	&20.90	&&11900	&20.38	&&69	&175.84	    &&3953	&298.49\\
&131p654	&654	&3	&&191	&24.71	&&9656	&24.83	&&98	&4542.93	&&4613	&4313.95\\
&132d657	&657	&3	&&202	&20.29	&&10174	&19.81	&&74	&12546.50	&&3811	&19328.10\\
&145u724	&724	&3	&&339	&29.25	&&16398	&29.14	&&78	&406.95	    &&4600	&972.32\\
&157rat783	&783	&3	&&333	&29.05	&&17023	&29.01	&&71	&521.87	    &&4150	&1191.50\\
&201pr1002	&1002	&3	&&304	&28.55	&&15040	&28.75	&&94	&2806.94	&&4871	&2599.29\\
&212u1060	&1060	&3	&&425	&37.49	&&21621	&37.16	&&103	&1608.19	&&5063	&1505.21\\
&217vm1084	&1084	&3	&&549	&53.80	&&27745	&54.02	&&145	&298.21    	&&7692	&326.57\\

\textbf{Avg.} &&&&& \textbf{314.67} &\textbf{29.76} &&\textbf{15841.22} & \textbf{29.66}        &  & \textbf{88.89}   & \textbf{2596.99}      &  & \textbf{4732.78}   & \textbf{3443.13} \\

&115rat575	&575	&4	&&310	&26.38	&&15270	&26.21	&&78	&723.27	    &&3893	&745.72\\
&115u574	&574	&4	&&269	&23.16	&&13491	&23.70	&&80	&518.63	    &&4454	&456.93\\
&131p654	&654	&4	&&205	&23.89	&&10219	&24.21	&&98	&19258.90   &&5195	&13670.80\\
&132d657	&657	&4	&&217	&20.42	&&10820	&20.48	&&74	&45989.70   &&4096	&47017.90\\
&145u724	&724	&4	&&419	&34.72	&&20866	&34.75	&&104	&1247.04	&&4600	&1379.36\\
&157rat783	&783	&4	&&441	&37.77	&&22289	&39.01	&&90	&1388.36	&&4334	&1011.20\\
&201pr1002	&1002	&4	&&465	&38.90	&&22246	&41.91	&&107	&1554.15	&&6165	&1189.69\\
&212u1060	&1060	&4	&&515	&45.00	&&26096	&44.83	&&111	&2003.73	&&5542	&3053.60\\
&217vm1084	&1084	&4	&&673	&58.80	&&33434	&59.51	&&157	&261.48	    &&6723	&228.01\\

\textbf{Avg.} &&&&& \textbf{390.44} &\textbf{34.34} &&\textbf{19414.56} & \textbf{34.95}        &  & \textbf{99.89}   & \textbf{8105.03}      &  & \textbf{5000.22}   & \textbf{7639.24} \\
\hline
\end{longtable}

In Table \ref{tab: Table 2}, we observe that the average time taken by the VNS to solve Set $1$ is $4.97$ seconds in the case of the rule $g_1$ and $4.98$ seconds in the case of the rule $g_2$. In the same case, GA takes $9.54$ and $10.57$ seconds, respectively. For Set $1$, the average profit is also better when we use VNS. In the case of the rule $g_1$, VNS earned $5.08\%$ more profit than GA using two travelers, but if we use three travelers, then VNS earned $4.48\%$ more profit than GA. In the case of the rule $g_2$, VNS earned $3.21\%$ more profit than GA using two travelers. While using three travelers, VNS earned $3.43\%$ more profit than GA.  

In the case of Set $2$ of Table \ref{tab: Table 2}, VNS consistently outperforms GA regarding average profit. Specifically, with two travelers, VNS achieves an average profit that is $63.10\%$ more using the rule $g_1$ and $73.33\%$ more using the rule $g_2$ compared to GA. Similarly, VNS outperforms GA with three travelers by $82.64\%$ using the rule $g_1$ and by $82.33\%$ using the rule $g_2$. 

In Set $3$ of Table \ref{tab: Table 2}, we use up to four travelers to observe the profit increment on large instances up to $1084$ nodes. There, we find a similar pattern after observing the results of VNS and GA. Specifically, using the rule $g_1$, VNS achieves $185.50\%$ more profit with two travelers, $253.999\%$ more profit with three travelers, and $290.43\%$ more profit with four travelers compared to GA. Similarly, with the rule $g_1$, VNS earned $170.23\%$ more profit with two travelers, $360.08\%$ more profit with three travelers, and $288.27\%$ more profit with four travelers than GA.

In Table \ref{tab: Table 2}, if we observe the profit gain collected by VNS using the different number of travelers for Set $1$, it collects $12.72\%$ more profit using three travelers instead of two using the rule $g_1$, and $10.67\%$ more profit using three travelers instead of two travelers using the rule $g_2$. A Similar pattern is observed in Set $2$ of Table \ref{tab: Table 2}; VNS collects $19.31\%$ more profit using three travelers instead of two travelers using the rule $g_1$ and $19.62\%$ more profit while using rule $g_2$ if three travelers are used instead of two travelers. However, in the case of large instances in Set $3$ using rule $g_1$, it earned $35.70\%$ more profit if three travelers are used instead of two travelers and $24.08\%$ more profit using four travelers instead of three travelers. In the case of the rule $g_2$, VNS collected $38.08\%$ more profit using three travelers instead of two travelers, and $22.56\%$ more profit using four travelers instead of three travelers.

\section{Conclusion}
\label{section:6}

This paper studies the single Depot multiple Set Orienteering Problem (sDmSOP), a multi-traveler variant of the Set Orienteering Problem (SOP), with a single starting and ending point. To solve this problem, we propose GA and VNS meta-heuristics. The results achieved on the benchmark instances for the sDmSOP show the effectiveness of our methods. Both methods show effective results on small instances. However, as the number of nodes increases, VNS consistently outperforms the GA method in terms of average profit earned from clusters and takes less computational time to solve the sDmSOP instances. We also observe that if a small value of $w$ is taken and the number of travelers is increased to solve the sDmSOP, it gives significantly better profit. If we take the larger value of $w$, the budget exceeds the solution of GTSP; in that case, only one traveler is required to visit all the clusters.

\bibliographystyle{apalike}

\bibliography{main}

\end{document}